# POSITIVE ASSOCIATION IN THE FRACTIONAL FUZZY POTTS MODEL[1]

By Jeff Kahn and Nicholas Weininger

*Rutgers University*

A *fractional fuzzy Potts measure* is a probability distribution on spin configurations of a finite graph $G$ obtained in two steps: first a subgraph of $G$ is chosen according to a random cluster measure $\phi_{p,q}$, and then a spin ($\pm 1$) is chosen independently for each component of the subgraph and assigned to all vertices of that component. We show that whenever $q \geq 1$, such a measure is positively associated, meaning that any two increasing events are positively correlated. This generalizes earlier results of Häggström [*Ann. Appl. Probab.* **9** (1999) 1149–1159] and Häggström and Schramm [*Stochastic Process. Appl.* **96** (2001) 213–242].

**1. Introduction.** We work with a finite graph $G = (V, E)$. A *random cluster measure* $\phi$ with parameter $q > 0$ is a probability measure on $\{0,1\}^E$ given by, for some weights $p_e \in (0,1)$ ($e \in E$),

$$\phi(\eta) \propto \prod_{e \in E} p_e^{\eta_e}(1-p_e)^{1-\eta_e} q^{k(\eta)}$$

where we interpret $\eta_e = 1$ as "$e$ is open," $k(\eta)$ is the number of connected components ("clusters") of the set of open edges specified by $\eta$, and, as usual, "$\propto$" means "proportional to." [One usually writes $\phi = \phi_{p,q}$ where $p = (p_e : e \in E)$.] For background on random cluster measures see, for example, Grimmett's survey [4].

A *fractional fuzzy Potts measure* is a probability measure $\nu = \nu_{\phi,\alpha} = \nu_{\phi,\alpha}^G$ on $\{\pm 1\}^V$ obtained by the following two-step process:

(i) Choose a random subgraph $\eta$ of $G$ with distribution $\phi = \phi_{p,q}$ (for some $p, q$).

Received July 2006; revised August 2006.

[1]Supported in part by NSF Grant DMS-02-00856.

*AMS 2000 subject classifications.* Primary 60C05; secondary 05D40.

*Key words and phrases.* Positive association, random cluster model, fuzzy Potts model.

(ii) Independently choose "spins" (1 or $-1$) for the clusters of $\eta$, where each spin is 1 with probability $\alpha \in (0, 1)$, and assign the spin of a cluster to each of its vertices.

Häggström [5] gives some motivation for considering such measures.

Recall that two events $A, B$ in some probability space are *positively correlated*—hereafter denoted by $A \uparrow B$—if $\Pr(AB) \geq \Pr(A) \Pr(B)$. The joint distribution of random variables $X_1, \ldots, X_n$ is said to be *positively associated* if any two events both increasing in the $X_i$'s are positively correlated. (This is easily seen to be equivalent to the property that for any two increasing functions $f, g$ of the $X_i$'s one has $\mathbf{E}fg \geq \mathbf{E}f\mathbf{E}g$.)

In this paper we show:

THEOREM 1. *Any fractional fuzzy Potts measure with $q \geq 1$ is positively associated.*

Some special cases were known earlier. In [5] Häggström proved Theorem 1 in case each of $\alpha q, (1 - \alpha)q$ is at least 1; he actually proved that in this case $\nu$ satisfies the stronger "positive lattice condition" [(1) below], and it is not hard to see that for this conclusion his conditions on $q, \alpha$ are also necessary. Later, Häggström and Schramm [6] proved Theorem 1 in the case $q = 1$, where the underlying $\phi$ is an ordinary bond percolation measure. We are now (weakly) conjecturing that Theorem 1 holds even when $q < 1$; see Section 3 for a little more on this.

(That (1)—also called the "FKG lattice condition"—implies positive association is the content of the celebrated "FKG inequality" of Fortuin, Kasteleyn and Ginibre [2]. This has been by far the most useful tool for proving positive association, in part because when (1) does hold, it is usually relatively easily seen to do so (though not always; see, e.g., [5]). Of course this means that, at least from a methodological standpoint, the most interesting positive association questions are those to which the FKG inequality does not apply.)

The proof of Theorem 1 is given in Section 2, and some additional remarks are included in Section 3.

**2. Proof of Theorem 1.** The proof of Theorem 1 is based on two (different) applications of the following, presumably standard observation, whose proof we omit.

LEMMA 1. *Suppose the events $A, B, C$ in some probability space satisfy:*

(i) *each of $A, B$ is positively correlated with $C$, and*
(ii) *$A$ and $B$ are conditionally positively correlated given either $C$ or $\overline{C}$.*

*Then $A$ and $B$ are positively correlated.*

We also need two well-known properties of random cluster measures (see, e.g., [4]). First, for any $e \in E$, the conditional measures $\phi_{p,q}(\cdot | \eta_e = 1)$ and $\phi_{p,q}(\cdot | \eta_e = 0)$ are random cluster measures with the same $q$ (and same $p_f$'s for $f \neq e$) on the graphs $G/e$, $G - e$, respectively. Second, for $q \geq 1$, $\phi_{p,q}$ satisfies the "positive lattice condition,"

$$(1) \qquad \forall \eta, \tau \in \{0,1\}^E \qquad \phi(\eta)\phi(\tau) \leq \phi(\eta \wedge \tau)\phi(\eta \vee \tau),$$

where $\wedge, \vee$ denote meet and join in the product order on $\{0,1\}^E$. We will use an equivalent version: for any $F \subseteq E$, $e \in E \setminus F$, and $\psi, \xi \in \{0,1\}^F$ with $\xi \leq \psi$,

$$(2) \qquad \phi(\eta_e = 1 \,|\, \eta \equiv \xi \text{ on } F) \leq \phi(\eta_e = 1 | \eta \equiv \psi \text{ on } F).$$

We now assume that $\alpha$, the random cluster measure $\phi = \phi_{p,q}$ $(q \geq 1)$ with associated edge configuration $\eta$, and the corresponding fractional fuzzy Potts measure $\nu$ are as described at the beginning of the paper, and write $\sigma = (\sigma_v : v \in V)$ for the random spin configuration produced by $\nu$. We must show $A \uparrow B$ for any increasing events $A, B$ determined by the $\sigma_v$'s.

This will follow easily from:

LEMMA 2.  *For any increasing event $C$ determined by the $\sigma_v$'s, $x \in V$, and $e \in E$ containing $x$, $C$ is conditionally positively correlated with the event $\{\eta_e = 1\}$ given $\{\sigma_x = 1\}$.*
(Note this is not true if $e$ does not contain $x$.)

PROOF.  It suffices to construct a coupled pair of random configurations $(\psi, \xi)$ of $E(G)$ such that:

(i)  $\psi$ has marginal distribution $\phi(\cdot | \eta_e = 1)$,
(ii)  $\xi$ has marginal distribution $\phi(\cdot | \eta_e = 0)$,
(iii)  $\Pr(C \,|\, \psi) \geq \Pr(C \,|\, \xi)$, these probabilities taken w.r.t. $\nu(\cdot | \sigma_x = 1)$.

We construct the coupling one edge at a time. Let $X_f$, $f \in E$, be independent (real-valued) random variables, each distributed uniformly on $[0,1]$. Set $e_0 = e$, $\psi_e = 1$ and $\xi_e = 0$, and repeat for $i = 1, \ldots, |E| - 2$:

Given $e_0, \ldots, e_i$ and the (random) restrictions, say $\psi_i, \xi_i$, of $\psi, \xi$ to $\{e_0, \ldots, e_i\}$, let $e_{i+1}$ be an edge of $G \setminus \{e_0, \ldots, e_i\}$ incident with the component of $x$ in $\psi_i$. If there is no such edge, let $e_{i+1}$ be any remaining edge. In either case, set $\psi_{i+1}(e_{i+1}) = 1$ iff $X_{e_{i+1}} \leq \phi(\eta(e_{i+1}) = 1 | \psi_i)$ and $\xi_{i+1}(e_{i+1}) = 1$ iff $X_{e_{i+1}} \leq \phi(\eta(e_{i+1}) = 1 | \xi_i)$, where the conditioning events are $\{\eta \equiv \psi_i \text{ on } \{e_0, \ldots, e_i\}\}$ (and similarly for $\xi_i$) and we use $\eta(f)$ in place of $\eta_f$, and so on.

Then $\psi, \xi$ clearly have the right marginals. Furthermore, if $\psi_i \geq \xi_i$, then $\psi_{i+1} \geq \xi_{i+1}$, since (2) implies that then

$$\phi(\eta(e_{i+1}) = 1 \,|\, \psi_i) \geq \phi(\eta(e_{i+1}) = 1 \,|\, \xi_i).$$

So indeed $\psi \geq \xi$. Finally, observe that, writing $S$ for the vertex set of the component of $x$ in $\psi$, $\psi$ and $\xi$ agree on all edges not contained in $S$; for, conditioned on the absence of (open) edges between $S$ and $\overline{S}$, the restriction of our random cluster measure to (edges contained in) $\overline{S}$ is independent of its restriction to $S$.

Thus every component of $\xi$ is either (a) identical to a component of $\psi$ or (b) contained in the $x$-component of $\psi$. But then we can couple $\nu(\cdot|\sigma_x = 1, \psi)$ and $\nu(\cdot|\sigma_x = 1, \xi)$ by choosing the same random spins for all components of $\psi$ other than the $x$-component; and this shows $\Pr(C\,|\,\psi) \geq \Pr(C\,|\,\xi)$ as desired. $\square$

PROOF OF THEOREM 1.  We proceed by induction on $|V| + |E|$, the case $|V| = 1$ being trivial. Fix some $x \in V$.

Observe first that $\{\sigma_x = 1\}$ is positively correlated with each of $A, B$ since for any $\eta$, $\Pr(A\,|\eta, \sigma_x = 1) \geq \Pr(A\,|\eta)$ and $\Pr(\eta\,|\sigma_x = 1) = \Pr(\eta) = \phi_{p,q}(\eta)$. So by Lemma 1 it suffices to show $A \uparrow B$ with respect to each of the conditional measures $\nu(\cdot|\sigma_x = 1)$ and $\nu(\cdot|\sigma_x = -1)$. Actually we only need to show this for $\{\sigma_x = 1\}$; symmetry then implies $\overline{A} \uparrow \overline{B}$ (which is equivalent to $A \uparrow B$) given $\{\sigma_x = -1\}$.

If $x$ is isolated, then positive correlation of $A$ and $B$ given $\{\sigma_x = 1\}$ is the same as positive correlation of $A' := \{\sigma|_{V \setminus \{x\}} : \sigma \in A, \sigma_x = 1\}$ and the analogously defined $B'$ under the measure $\nu_{\phi,\alpha}^{G-x}$ [note $\phi$ still makes sense here since $E(G - x) = E(G)$], so holds by inductive hypothesis.

Otherwise, let $e$ be an edge containing $x$. It follows from Lemma 2 that each of $A, B$ is positively correlated with $\{\eta_e = 1\}$ given $\{\sigma_x = 1\}$. And by inductive hypothesis we have $A \uparrow B$ with respect to either of $\nu(\cdot|\sigma_x = 1, \eta_e = 1), \nu(\cdot|\sigma_x = 1, \eta_e = 0)$, since $\phi(\cdot|\eta_e = 1), \phi(\cdot|\eta_e = 0)$ are random cluster measures on smaller graphs with $q \geq 1$. So Lemma 1 gives positive correlation of $A, B$ w.r.t. $\nu(\cdot|\sigma_x = 1)$, as desired.  $\square$

## 3. Questions for further investigation.

**A.** Note that in the preceding proof, we used only two properties of $\phi_{p,q}$:

  (i)  the positive lattice condition;
  (ii)  the fact that if $S \subseteq V(G)$, and $F$ denotes the event that there are no (open) edges joining $S$ and $\overline{S}$, then the restrictions of $\phi_{p,q}(\cdot\,|\,F)$ to $S$ and $\overline{S}$ are independent.

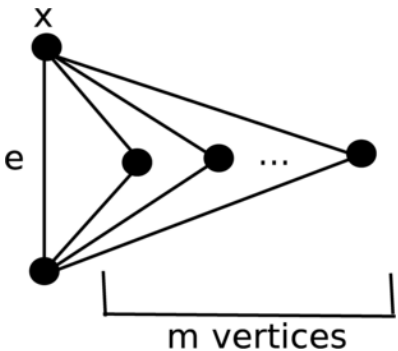

Fig. 1. *Lemma 2 may not hold when $q < 1$.*

So we can replace $\phi_{p,q}$ by any other measure $\phi$ having these properties, and the resultant measure $\nu_{\phi,\alpha}$ will be positively associated. It is easy to construct toy measures $\phi$ that have these properties but are not random cluster measures, but it would be interesting to see if any such measures arise naturally in other probabilistic contexts. It is also easy to show by small example that neither property alone suffices to give positive association of $\nu$.

**B.** As mentioned above, we think Theorem 1 may hold even when $q < 1$. This would be quite striking, since the correlation properties of $\phi_{p,q}$ for $q < 1$ are usually completely different from—and much more mysterious than—those for $q \geq 1$. For instance, except in trivial cases (graphs without cycles), random cluster measures with $q < 1$ fail to have the positive lattice condition or even positive association. In fact, negative association properties are thought to hold; see [4] for details.

More significantly for the present discussion, when $q < 1$ we lose Lemma 2. To see this, first observe that (as is well known) if we let $p_e = p = q$ for each edge $e$ and let $q \to 0$, then $\phi_{p,q}$ approaches uniform measure on spanning forests of $G$. So it is enough to show that Lemma 2 can fail when $\nu = \nu_{\phi,\alpha}$ corresponds to a uniform forest measure $\phi$.

Now consider $\nu_{\phi,\alpha}$ on the graph shown in Figure 1. Let $\alpha$ approach zero and let $A$ be the event that all vertices are colored 1. Then $\Pr(A \,|\, \sigma_x = 1)$ approaches the probability that the underlying subgraph $\eta$ is connected (since, on $\{\sigma_x = 1\}$, $A$ always occurs if $\eta$ is connected and has probability at most $\alpha$ otherwise). But a simple calculation shows that if $m > 6$, the events $\{\eta_e = 1\}$ and $\{\eta$ is connected$\}$ are strictly negatively correlated. (Related properties of this graph have probably been rediscovered several times, but as far as we know were first observed by Dilworth and Greene [1].)

**C.** There is a natural, more general context for the present discussion. For any finite set $V$ and any probability measure $\psi$ on the set of unordered partitions of $V$, we may construct a measure $\nu_{\psi,\alpha}$ on $\{\pm 1\}^V$ by choosing a partition $\pi$ according to $\psi$ and then choosing a spin for each "cluster" (i.e., block) of $\pi$ as before. We might then ask what other conditions on $\psi$ would suffice for $\nu$ to be positively associated.

One may also study analogous processes with more colors. Let $V, \psi, \pi$ be as in the preceding paragraph and assign colors from some finite set $\Omega$ to the blocks of $\pi$, these colors chosen independently, each according to some fixed distribution $\beta$ on $\Omega$. We would like to use the term "divide and color" for such procedures, following [6] where it was used for the case in which the blocks of $\pi$ are the clusters of some i.i.d. bond percolation. Another special case is the $q$-state *Potts model* (for which see, e.g., [3]), obtained when the blocks are the clusters of $\eta$, which is the output of some random cluster measure with $q \geq 2$ an integer, and $\beta$ is uniform on $\{0, \ldots, q-1\}$. Also of this type are the processes of [7], Section 3, though the concerns there are rather different than ours.

**Acknowledgments.** We thank Mike Neiman for pointing out an error in an earlier version of the paper, and the referee for drawing our attention to [7].

DEPARTMENT OF MATHEMATICS
RUTGERS UNIVERSITY
PISCATAWAY, NEW JERSEY 08854
USA
E-MAIL: jkahn@math.rutgers.edu
        nweining@math.rutgers.edu